\documentclass[11pt]{article}
\usepackage{amsmath}
\usepackage{amsfonts}
\newtheorem{thm}{Theorem}[section]
\newtheorem{lem}[thm]{Lemma}

\newtheorem{cor}[thm]{Corollary}
\newtheorem{defn}[thm]{Definition}
\newenvironment{defn-new}{\begin{defn} \em}{\end{defn}}
\newtheorem{rem}[thm]{Remark}
\newenvironment{rem-new}{\begin{rem} \em}{\end{rem}}
\newtheorem{ex}[thm]{Example}
\newenvironment{ex-new}{\begin{ex} \em}{\end{ex}}

\newenvironment{notation-new}{\begin{rem} \em}{\end{rem}}

\newenvironment{agr-new}{\begin{rem} \em}{\end{rem}}

\makeatletter \@addtoreset{equation}{section} \makeatother

\makeatletter \@addtoreset{figure}{section} \makeatother

\begin{document}

\begin{center}
{\bf {\Large Curvature inequalities for anti-invariant  submersion from quaternionic space forms}}\\[0pt]

{\bf Kirti Gupta}\footnote{
Department of Mathematics \& Statistics, Dr. Harisingh Gour Vishwavidyalaya,
Sagar-470 003, M.P. INDIA \newline
Email: guptakirti905@gmail.com}, {\bf Punam Gupta}\footnote{
School of Mathematics, Devi Ahilya Vishwavidyalaya, Indore-452 001, M.P. 
INDIA\newline
Email: punam2101@gmail.com} and {\bf R.K. Gangele}\footnote{
Department of Mathematics \& Statistics, Dr. Harisingh Gour Vishwavidyalaya,
Sagar-470 003, M.P. INDIA \newline
Email: rkgangele23@gmail.com} 
\end{center}

\noindent {\bf Abstract:}  This paper focuses on deriving several curvature inequalities involving the Ricci and scalar curvatures of the horizontal and vertical distributions in anti-invariant Riemannian submersions from quaternionic space forms onto Riemannian manifolds. In addition, a Ricci curvature inequality for anti-invariant Riemannian submersions is established. The equality cases for all derived inequalities are also examined.


\section{Introduction}
In 1993, B.-Y. Chen established foundational relationships between the extrinsic and intrinsic invariants of submanifolds in real space forms \cite{Chen}. Later, in 1999, he derived a sharp inequality connecting the Ricci curvature and the squared mean curvature for submanifolds \cite{Chen99}. These findings have since inspired a broad spectrum of research, leading to the development of numerous results in various ambient geometries. For a comprehensive study, see \cite{Chenbook}. Several authors \cite{Aydin, Ayi,  Naghi, Sahin16, MMT} have further explored Chen-like inequalities.

As noted in \cite{Falcitelli-Ianus and Pastore}, a central theme in Riemannian geometry is the comparison of geometric properties of different types of maps between Riemannian manifolds. In this context, the concept of Riemannian submersion, introduced by B. O’Neill in 1966 \cite{O'Neill}, plays a pivotal role. O'Neill formulated key equations governing Riemannian submersions, which have found applications in diverse fields such as supergravity, statistical machine learning, analysis on manifolds, and image processing. The geometric features of Riemannian submersions have been further investigated by Besse \cite{Besse}, Escobales Jr. \cite{Eco}, and many others.

In a related development, Watson \cite{Watson} introduced and studied almost Hermitian submersions between almost complex manifolds, showing that both the horizontal and vertical distributions remain invariant under the almost complex structure of the total space. Later, in 2010, \c Sahin \cite{Sahin-10} explored anti-invariant Riemannian submersions from almost Hermitian manifolds to Riemannian manifolds, demonstrating that the fibers are anti-invariant under the almost complex structure. Consequently, the horizontal distribution is not preserved by this structure.

In this manuscript, we investigate inequalities involving the Ricci and scalar curvatures on the vertical and horizontal distributions of anti-invariant Riemannian submersions from quaternionic space forms. Section 2 presents the foundational geometric properties of Riemannian submersions and quaternionic K\"ahler manifolds. In Section 3, we derive new inequalities for the Ricci and scalar curvatures associated with the vertical and horizontal distributions, and we analyze the equality cases in detail.
\section{Preliminaries}
This section revisits fundamental definitions and presents established results related to Riemannian submersions and quaternionic K\"ahler manifolds, which will be utilized in the paper.  

\subsection{Submersion}

We will now present a set of definitions and results for future reference.

\begin{defn-new}
{\rm \cite{Baird}} A $C^{\infty }$-map $\pi:\left(
M,g_{M}\right)  \rightarrow \left( N,g_{N}\right) $  between Riemannian  manifolds $\left( M,g_{M}\right) $ and $\left( N,g_{N}\right) $ is said to be  $C^{\infty }$-submersion if
\begin{description}
\item[(i)] $\pi$ is surjective. 
\item[(ii)] The differential $\left(\pi_*\right)_p$ has maximal rank for any $
p \in M$.\newline
\end{description}
\end{defn-new}
\begin{defn-new}
{\rm {\cite{Falcitelli-Ianus and Pastore, O'Neill}}}
Let $(M,g_{M})$ and $
(N,g_{N})$ be Riemannian manifolds with $\dim (M)=m>n=\dim (N)$. A smooth surjective map $\pi :M\rightarrow N$ is said to be Riemannian submersion if $\pi$ satisfies the
following axioms:
\begin{description}
\item[(i)] $\pi$ has maximal rank.
\item[(ii)] The differential $\pi_{\ast }|_{(\operatorname{ker} \pi
_{\ast })^{\perp }}$ is a linear isometry.
\end{description}
\end{defn-new}

\noindent{\bf Remark:} For each $q\in N$, $\pi ^{-1}(q)$ is an $(m-n)$-dimensional submanifold of $%
M $. The submanifolds $\pi ^{-1}(q)$, $q\in N$, are called fibers. A vector
field on $M$ is called vertical if it is always tangent to the fibers. A vector
field on $M$ is called horizontal if it is always orthogonal to the fibers. A
vector field $X$ on $M$ is called basic if $X$ is horizontal and $\pi $%
-related to a vector field $X^{\prime }$ on $N$, that is, $\pi _{\ast
}X_{p}=X_{\pi _{\ast }(p)}^{\prime }$ for all $p\in M.$ We denote the
projection morphisms on the distributions $\operatorname{ker}\pi _{\ast }$ and $(\operatorname{ker} \pi
_{\ast })^{\perp }$ by ${\cal V}$ and ${\cal H}$, respectively. The sections
of ${\cal V}$ and ${\cal H}$ are called vertical vector fields and
horizontal vector fields, respectively. So
\[
{\cal V}_{p}=T_{p}\left( \pi ^{-1}(q)\right) ,\qquad {\cal H}%
_{p}=T_{p}\left( \pi ^{-1}(q)\right) ^{\perp }.
\]

\medskip

\noindent Let $\pi:\left( M,g_{M}\right) \rightarrow \left( N,g_{N}\right) $ 
be a Riemannian submersion. Given any vector field $U\in \Gamma (TM)$, 
we have  
\[
U={\cal V}U+{\cal H}U,  
\]
where ${\cal V}U\in \Gamma \left( {\operatorname{ker}}\pi_{\ast }\right)$ and
${\cal H}U\in \Gamma \left( \left( {\operatorname{ker}}\pi_{\ast }\right) ^{\perp }\right) $.

The second fundamental tensors of all fibers $\pi ^{-1}(q),\ q\in N$ give
rise to configuration tensors ${\mathcal{T}}$ and ${\mathcal{A}}$ in $M$ defined by O'Neill 
\cite{O'Neill}.  
\begin{eqnarray}
{\cal A}_{E}F &=&{\cal H}\nabla _{{\cal H}E}{\cal V}F+{\cal V}\nabla _{{\cal %
H}E}{\cal H}F, \\
{\cal T}_{E}F &=&{\cal H}\nabla _{{\cal V}E}{\cal V}F+{\cal V}\nabla _{{\cal %
V}E}{\cal H}F
\end{eqnarray}
for $E,F\in \Gamma (TM)$, where $\nabla $ is the Levi-Civita connection of $
g_{M}$.

From the above equations, we have  
\begin{equation}
\nabla _{X}Y={\cal T}_{X}Y+\widehat{\nabla }_{X}Y,
\end{equation}
\vspace{-0.7cm}
\begin{equation}
\nabla _{X}V={\cal H}\nabla _{X}V+{\cal T}_{X}V ,  \label{eq-aa}
\end{equation}
\begin{equation}
\nabla _{V}X={\cal A}_{V}X+{\cal V}\nabla _{V}X,
\end{equation}
\begin{equation}
\nabla _{V}W={\cal H}\nabla _{V}W+{\cal A}_{V}W,
\end{equation}
for $X,Y\in \Gamma  \left( {\operatorname{ker}}\pi_{\ast }\right)$ and $V,W\in
\Gamma \left( \left( {\operatorname{ker}}\pi_{\ast }\right) ^{\perp }\right) $, where $\widehat{
\nabla}_{X}Y={\cal V}\nabla _{X}Y$. If $V$ is basic, then ${\mathcal{A}}_{V}X={\cal H} 
\nabla _{X}V.$

\noindent It is easy to find that for $p\in M,$ $X\in {\cal V}_{p}$ and $U\in {\cal H}
_{p}$, the linear operators  
\[
{\mathcal{A}}_{U},{\mathcal{T}}_{X}:T_{p}M\rightarrow T_{p}M  
\]
are skew-symmetric. 
\newline The tensor fields ${\mathcal{T}}$ and ${\mathcal{V}}$ satisfy
\begin{align}
    {\mathcal{T}}_{U}W={\mathcal{T}}_{W}U,\\
    {\mathcal{A}}_{X}Y=-{\mathcal{A}}_{Y}X
\end{align}
for any $U,W \in {\mathcal{V}}(M)$ and $X,Y \in {\mathcal{H}}(M).$
\begin{lem} \rm{\cite{Gulbahar}}
    Let $ (M, g_{M})$  and $(N,g_{N})$ be Riemannian manifolds admitting a Riemannian submersion $\pi:M \rightarrow N.$ Let  $\left\{U_1, ..., U_r
, X_1, ..., X_\ell \right \}$  be an orthonormal
frame of $T_{p}M$  such that ${\mathcal{V}} = span \left \{U_1, ..., U_r  \right \}, {\mathcal{H}} = span \left \{X_1, ..., X_\ell \right \}$. Then at any point $p \in M,$ we have
    \begin{equation} \label{eq-3.15}
\sum_{s=1}^{\ell} \sum_{i,j=1}^{r}{({\mathcal{T}}_{ij}^{s})}^{2}= \frac{1}{2}r^{2} \lVert H \rVert ^{2} + ({\mathcal{T}}_{11}^{s}-{\mathcal{T}}_{22}^{s}- ... - {\mathcal{T}}_{rr}^{s})^{2} + 2\sum_{s=1}^{\ell} \sum_{j=2}^{r} ({\mathcal{T}}_{1j}^{s})^{2}-2\sum_{s=1}^{\ell} \sum_{2\leq i<j \leq r} ({\mathcal{T}}_{ii}^{s}{\mathcal{T}}_{jj}^{s}-({\mathcal{T}}_{ij}^{s})^{2}).
\end{equation}
\end{lem}

\subsection{Quaternionic K\"ahler Manifold}

Let $M$ be a $4m$-dimensional $C^{\infty }$-manifold and let $E$ be a rank 3
subbundle of End $(TM)$ such that for any point $p\in M$ with a neighborhood
$U$, there exists a local basis $\left\{ J_{1},J_{2},J_{3}\right\} $ of 
sections of $E$ on $U$ satisfying  
\[
J_{\alpha}^{2}=-id,\quad J_{\alpha }J_{\alpha +1}=-J_{\alpha +1}J_{\alpha
}=J_{\alpha +2}  
\]
for all $\alpha \in \{1,2,3\}$,  where the indices are taken from $\{1,2,3\}$
modulo $3$. Then $E$ is said to  be an almost quaternionic structure on $M$
and $(M,E)$ an almost  quaternionic manifold \cite{Alekseevsky-Marchiafava}.

Moreover, let $g_M$ be a Riemannian metric on $M$ defined by  
\[
g_M\left( J_{\alpha }X,J_{\alpha }Y\right) =g_M(X,Y)  
\]
for all vector fields $X,Y\in \Gamma (TM)$, where the indices are taken from
$\{1,2,3\}$ modulo $3$. Then $(M,E,g_M)$ is said to be an almost quaternionic Hermitian manifold \cite%
{Ianus-Mazzocco-Vilcu} and the basis $\left\{ J_{1},J_{2},J_{3}\right\} $ is
said to be a quaternionic Hermitian basis.

An almost quaternionic Hermitian manifold $(M,E,g_M)$ is said to be a 
quaternionic K\"{a}hler manifold \cite{Ianus-Mazzocco-Vilcu} if there exist 
locally defined $1$-forms $\omega _{1},\omega _{2},\omega _{3}$ such that 
for $\alpha \in \{1,2,3\}$ modulo $3$  
\begin{equation}
\nabla _{X}J_{\alpha }=\omega _{\alpha +2}(X)J_{\alpha+1}-\omega _{\alpha
+1}(X)J_{\alpha +2}  \label{eq-qkm}
\end{equation}
for $X\in \Gamma (TM)$, where the indices are taken from $\{1,2,3\}$ modulo $
3$.

If there exists a global parallel quaternionic Hermitian basis $\left\{ 
J_{1},J_{2},J_{3}\right\} $ of sections of $E$ on $M$ (i.e. $\nabla J_{\alpha}=0$
for $\alpha \in \{1,2,3\}$, where $\nabla $ is the Levi-Civita connection of
the metric $g_M$, then $(M,E,g_M)$ is said to be a hyperk\"ahler manifold \cite%
{Besse}.
\begin{defn}
We call a Riemannian submersion $\pi:\left( M,E,g_{M}\right) \rightarrow 
\left( N,g_{N}\right) $ from quaternionic K\"ahler 
manifold $(M,E,g_{M})$  onto a Riemannian manifold $\left( N,g_{N}\right) $ an anti-invariant submersion {\rm \cite{Park-17}}
 if $J_\alpha\left( \operatorname{ker}\pi_{\ast }\right) \subset \left( {\ \operatorname{ker}} 
\pi_{\ast }\right) ^{\perp }$. 
\end{defn}
\noindent Let $\pi$ be an anti-invariant submersion  from quaternionic K\"ahler 
manifold $(M,E,g_{M})$  onto a Riemannian manifold $\left( N,g_{N}\right) $%
. Given a point $p\in M$ with a neighborhood $U$, we have an  
anti-invariant basis $
\{J_{1},J_{2},J_{3}\}$ of sections of $E$ on $U$. Given $X\in \Gamma  \left(
\left( \operatorname{ker}\pi_{\ast }\right) ^{\perp }\right) $, we have  
\begin{equation}\label{eq-7..}
J_\alpha X=B_{J_\alpha}X+C_{J_\alpha}X, \quad \alpha \in \{1,2,3\},
\end{equation}
where $B_{J_\alpha}X\in \Gamma \left( \operatorname{ker}\pi_{\ast }\right) $ and $C_{J_\alpha}X\in \Gamma
\left( \left( \operatorname{ker}\pi_{\ast }\right) ^{\perp }\right) $.

We denote the Riemannian curvature tensors of Riemannian manifolds $(M,g_{M}), (N,g_{N}),$ the vertical distribution ${\mathcal{V}}$ and the horizontal distribution ${\mathcal{H}},$ by $R, R', \hat{R}, R^{*}$, respectively.
Then the Gauss-Codazzi type equations \cite{O'Neill} are given by the following:
\begin{eqnarray}
    R(U,V,F,W) = \hat{R}(U,V,F,W) + g_{M}({\mathcal{T}}_{U}W,{\mathcal{T}}_{V}F) - g_{M}({\mathcal{T}}_{V}W,{\mathcal{T}}_{U}F),\label{eq-1} 
    \end{eqnarray}
    \begin{eqnarray}
    R(X,Y,Z,H) &=& R^{*}(X,Y,Z,H) - 2g_{M}({\mathcal{A}}_{X}Y,{\mathcal{A}}_{Z}H)\notag\\&& +g_{M}({\mathcal{A}}_{Y}Z,{\mathcal{A}}_{X}H) - g_{M}({\mathcal{A}}_{X}Z,{\mathcal{A}}_{Y}H), \label{eq-2} 
    \end{eqnarray}
    \begin{eqnarray}
    R(X,V,Y,W) &=& g_{M}((\nabla_{X}{\mathcal{T}})(V,W),Y) + g_{M}((\nabla_{V}{\mathcal{A}})(X,Y),W)\notag\\&&  - g_{M}({\mathcal{T}}_{V}X,{\mathcal{T}}_{W}Y) + g_{M}({\mathcal{A}}_{Y}W,{\mathcal{A}}_{X}V),\label{eq-3}
\end{eqnarray}
and
\begin{equation}
    \pi_{*}(R^{*}(X,Y)Z)=R'(\pi_{*}X,\pi_{*}Y)\pi_{*}Z,
\end{equation}
for all $U,V,F,W \in {\mathcal{V}}(M)$ and $X,Y,Z,H \in {\mathcal{H}}(M).$

\vspace{0.4cm}
\noindent Let $\{U_{1},\ldots,U_{r}\}$ be an orthonormal frame of the vertical distribution ${\mathcal{V}}$. The mean curvature vector field $H$ of any fibre of the Riemannian submersion $\pi$ is given by
\begin{equation}
    H=\frac{1}{r}N, \quad {\text \rm where }  \quad N=\sum_{j=1}^{r}{\mathcal{T}}_{U_{j}}U_{j}.
\end{equation}
The Riemannian submersion $\pi$ has totally geodesic fibres if ${\mathcal{T}}$ vanishes identically on $M$ and $\pi$ has totally umbilical fibres if
\[ {\mathcal{T}}_U V ={\mathcal T}(U,V)= g_{M}(U, V) H, \]
where $U, V \in {\mathcal{V}}(M)$ and $H$ is the mean curvature vector field of any fibre. 
The horizontal distribution $\mathcal{H}$ is integrable if  ${\mathcal{A}}$  vanishes identically. 

\vspace{0.2cm}
\noindent Let $\left( M,E,g_{M}\right) $ be a quaternionic K\"ahler manifold and let $X$ be a non-null vector on $M.$
Then the $4$-plane spanned by $\{X,J_{1}X,J_{2}X,J_{3}X\}$, denoted by $Q(X),$  is called a quaternionic
$4$-plane. Any $2$-plane in $Q(X)$  is called a quaternionic plane. The sectional curvature of
a quaternionic plane is called a quaternionic sectional curvature. A quaternionic K\"ahler
manifold is a quaternionic space form if its quaternionic sectional curvatures are equal
to a constant, say $c.$ A quaternionic K\"ahler manifold $ (M, E, g_{M})$  is a
quaternionic space form, denoted $M(c),$  if and only if its Riemannian curvature tensor $R$ on $M(c)$ is given by \rm{\cite{S Ishi}}
\vspace{-.1cm}
\begin{equation} \label{eq-2.9}
\begin{aligned}
    R(X, Y)Z = \frac{c}{4} \Bigg\{ & g_{M}(Y, Z)X - g_{M}(X, Z)Y \\
    & + \Bigg[ \sum_{\alpha = 1}^{3} \Big(
        g_{M}(J_{\alpha}Y, Z) J_{\alpha}X 
        - g_{M}(J_{\alpha}X, Z) J_{\alpha}Y 
        + 2g_{M}(J_{\alpha}Y, X) J_{\alpha}Z
    \Big) \Bigg] \Bigg\},
\end{aligned}
\end{equation}
for $X,Y,Z\in \Gamma (TM)$.

\subsection{Some useful results}
Now, we are stating some results that we will use in this work.
\begin{lem} \label{lem-14}{\rm \cite{Apostol}} Let $a_1, a_2, \ldots , a_n$ be $n$-real numbers ( $n>1$ ), then

$$
\frac{1}{n}\left(\sum_{i=1}^n a_i\right)^2 \leq \sum_{i=1}^n a_i^2
$$
with equality if and only if  $a_1=a_2=\ldots=a_n$.
\end{lem}
\begin{lem} \label{lem-11} Let $a$ and $b$ be non-negative real numbers, then

$$
\frac{a+b}{2} \geq \sqrt{a b}
$$
with equality if and only if $a=b$.
\end{lem}
\section{Ricci Inequalities for Quaternionic Space Form}

In this section, we establish some inequalities on the vertical and horizontal distributions for anti-invariant submersions from quaternionic space forms to Riemannian manifolds.\\

\noindent Let  $(M(c), g_{M})$ be a quaternionic space form, $ (N, g_{N})$ a Riemannian manifold and
 $\pi : M(c) \rightarrow N$  an anti-invariant submersion. Furthermore, let  $\left\{U_1, ..., U_r
, X_1, ..., X_\ell \right \}$  be an orthonormal
frame of $T_{p}M(c)$  such that ${\mathcal{V}} = span \left \{U_1, ..., U_r  \right \}, {\mathcal{H}} = span \left \{X_1, ..., X_\ell \right \}.$  Then using \rm(\ref{eq-1}), \rm(\ref{eq-2}) and \rm(\ref{eq-3}), we have
\begin{equation}
\begin{aligned}
    \widehat{R}(U,V,F,W) &=\frac{c}{4} \Bigg\{ g_{M}(V,F)g_{M}(U,W) - g_{M}(U,F)g_{M}(V,W) \\
    &+ \Bigg[ \sum_{\alpha = 1}^{3} \Big( g_{M}(J_{\alpha}V,F) g_{M}(J_{\alpha}U,W) 
        - g_{M}(J_{\alpha}U, F) g_{M}(J_{\alpha}V,W) \\
        &+ 2g_{M}(J_{\alpha}V, U) g_{M}(J_{\alpha}F,W)
    \Big) \Bigg] \Bigg\}- g_{M}({\mathcal{T}}_{U}W,{\mathcal{T}}_{V}F)+g_{M}({\mathcal{T}}_{V}W,{\mathcal{T}}_{U}F),
\end{aligned}
\label{3.1}
\end{equation}

\begin{equation}
\begin{aligned}
    R^{*}(X,Y,Z,H) &=\frac{c}{4} \Bigg\{ g_{M}(Y,Z)g_{M}(X,H) - g_{M}(X,Z)g_{M}(Y,H) \\
    &+ \Bigg[ \sum_{\alpha = 1}^{3} \Big( g_{M}(J_{\alpha}Y,Z) g_{M}(J_{\alpha}X,H) 
        - g_{M}(J_{\alpha}X, Z) g_{M}(J_{\alpha}Y,H) \\
        &+ 2g_{M}(J_{\alpha}Y, X) g_{M}(J_{\alpha}Z,H)
    \Big) \Bigg] \Bigg\}+ 2g_{M}({\mathcal{A}}_{X}Y,{\mathcal{A}}_{Z}H)-g_{M}({\mathcal{A}}_{Y}Z,{\mathcal{A}}_{X}H) \\
    &+g_{M}({\mathcal{A}}_{X}Z, {\mathcal{A}}_{Y}H).
\end{aligned}
\label{eq-3.2}
\end{equation}

\begin{thm}
Let $ \pi : M(c) \rightarrow N $  be an anti-invariant submersion from a quaternionic space form $ (M(c), g_M) $ onto
a Riemannian manifold $ (N, g_N).$ Then  
\begin{equation} \label{eq-3.3}
 \widehat{Ric}(U) \geq \frac{c}{4}(r-1)-r g_{M}\left( {\mathcal{T}}_{U}U,H\right), 
\end{equation} 
where $U$ is a unit vertical vector of $M(c)$.\newline
The equality case of {\rm (\ref{eq-3.3})} holds if and only if each fiber is totally geodesic.
\end{thm}

\noindent {\bf Proof:} Let $\{U_{1},U_{2}\ldots,U_{r}\}$ be a local orthonormal frame of the vertical distribution $ \operatorname{ker} \pi_{\ast}$ and $U$ a vertical vector. 
Since 
\[\widehat{Ric}(U)\equiv\widehat{Ric}(U,U) =\sum_{i = 1}^{r} \hat{R}(U,U_{i},U_{i},U).\]
Using equation \rm{(\ref{3.1})} we have 
$$
\begin{aligned}
    \widehat{Ric}(U) &=\frac{c}{4} \sum_{i = 1}^{r}\Bigg\{ g_{M}(U_{i},U_{i})g_{M}(U,U) - g_{M}(U,U_{i})g_{M}(U_{i},U) \\
    &+ \Bigg[ \sum_{\alpha = 1}^{3} \Big( g_{M}(J_{\alpha}U_{i},U_{i}) g_{M}(J_{\alpha}U,U) 
        - g_{M}(J_{\alpha}U, U_{i}) g_{M}(J_{\alpha}U_{i},U) \\
        &+ 2g_{M}(J_{\alpha}U_{i}, U) g_{M}(J_{\alpha}U_{i},U)
    \Big) \Bigg] \Bigg\}-\sum_{i = 1}^{r} g_{M}({\mathcal{T}}_{U}U,{\mathcal{T}}_{U_{i}}U_{i})+\sum_{i = 1}^{r}g_{M}({\mathcal{T}}_{U_{i}}U,{\mathcal{T}}_{U}U_{i}).
\end{aligned}
$$
Since $\pi$ is an anti-invariant submersion, therefore we get $g_{M}\left( J_{\alpha}U,U_{i}\right)=0 \quad \forall i, \alpha.$ \\
Now, we have 
\[ \widehat{Ric}(U)= \frac{c}{4}\sum_{i = 1}^{r}\left( 1-g_{M}(U,U_{i})^2\right)-\sum_{i = 1}^{r} g_{M}({\mathcal{T}}_{U}U,{\mathcal{T}}_{U_{i}}U_{i})+\sum_{i = 1}^{r}g_{M}({\mathcal{T}}_{U_{i}}U,{\mathcal{T}}_{U}U_{i}).\]
\begin{equation} \label{eq-6}
\widehat{Ric}(U)= \frac{c}{4}(r-1)-rg_{M}\left( {\mathcal{T}}_{U}U,H\right)+ \sum_{i = 1}^{r} \lVert{\mathcal{T}}_{U}U_{i} \rVert ^{2}. 
\end{equation}
\[\widehat{Ric}(U)\geq  \frac{c}{4}(r-1)-rg_{M}\left( {\mathcal{T}}_{U}U,H\right).\]
By using (\ref{eq-6}), we can say that the equality case of {\rm (\ref{eq-3.3})} holds if and only if each fiber is totally geodesic.
\begin{thm}
Let $ \pi : M(c) \rightarrow N $  be an anti-invariant submersion from a quaternionic space form $ (M(c), g_M) $ onto
a Riemannian manifold $ (N, g_N).$ Then
\begin{equation} \label{eq-3.4}
\hat{\tau} \geq \frac{c}{4}r(r-1) - r^{2} \lVert H\rVert ^{2}.
\end{equation}
The equality case of \rm (\ref{eq-3.4}) holds if and only if each fiber is totally geodesic.
\end{thm}
\noindent {\bf Proof:} 
We know that 
\[ \hat{\tau} = \sum_{j=1}^{r} \widehat{Ric}(U_{j},U_{j}).\]
Using \rm(\ref{eq-6}), we get
\begin{equation} \label{eq-8.} 
\hat{\tau} = \frac{c}{4}{r(r-1)}- r^{2} \lVert H \rVert ^{2}+ \sum_{i,j = 1}^{r} \lVert{\mathcal{T}}_{U_j}U_{i} \rVert ^{2} .\end{equation}
\[  \hat{\tau} \geq \frac{c}{4}{r(r-1)}- r^{2} \lVert H \rVert ^{2}.\]

\begin{thm}
Let $\pi: M(c) \rightarrow N $ be an anti-invariant submersion from a quaternionic space form $(M(c),g_{M})$ onto a Riemannian manifold $(N,g_{N}).$ Then
\begin{equation} \label{eq-7}
Ric^*(X) \leq \frac{c}{4} \Bigg\{ (\ell -1) 
    +\sum_{j=1}^{\ell} \sum_{\alpha = 1}^{3} 
        3g_{M}(C_{J_{\alpha}}X, X_{j})^2
      \Bigg\},\end{equation}
where $X$ is a unit horizontal vector of $M(c)$.\\
The equality case of \rm(\ref{eq-7}) holds if and only if ${\mathcal{H}}(M)$ is integrable.
\end{thm}
{\bf Proof:}  Let $\{X_{1},X_{2}\ldots,X_{\ell}\}$ be a local orthonormal frame of the horizontal distribution $ (\operatorname{ker} \pi_{\ast})^\perp$.
Since 
\[ {Ric}^{*}(X) =\sum_{j = 1}^{\ell} {R}^{*}(X,X_{j},X_{j},X),\]
using equation \rm{(\ref{eq-3.2})}, we have 
\begin{equation}
\begin{aligned} \label {eq-3.7}
    Ric^*(X) &=\sum_{j=1}^{_\ell}\frac{c}{4} \Bigg\{ g_{M}(X_{j},X_{j})g_{M}(X,X) - g_{M}(X,X_{j})g_{M}(X_{j},X) \\
    &+ \Bigg[ \sum_{\alpha = 1}^{3} \Big( g_{M}(J_{\alpha}X_{j},X_{j}) g_{M}(J_{\alpha}X,X) 
        - g_{M}(J_{\alpha}X, X_{j}) g_{M}(J_{\alpha}X_{j},X) \\
        &+ 2g_{M}(J_{\alpha}X_{j}, X) g_{M}(J_{\alpha}X_{j},X)
    \Big) \Bigg] \Bigg\}+ 2g_{M}({\mathcal{A}}_{X}X_{j},{\mathcal{A}}_{X_{j}}X)\\
    &-g_{M}({\mathcal{A}}_{X_{j}}X_{j},{\mathcal{A}}_{X}X) +g_{M}({\mathcal{A}}_{X}X_{j}, {\mathcal{A}}_{X_{j}}X).
\end{aligned}
\end{equation}
Using (\ref{eq-7..}), we have $$ 
\begin{aligned}
  Ric^*(X) &=\sum_{j=1}^{\ell}\frac{c}{4} \Bigg\{ (1-g_{M}(X,X_{j})^{2}) 
    + \sum_{\alpha = 1}^{3} 
        3g_{M}(C_{J_{\alpha}}X, X_{j})^2
      \Bigg\}\\&
    - \sum_{j=1}^{\ell}\left(3 g_{M}({\mathcal{A}}_{X}X_{j},{\mathcal{A}}_{X}X_{j})
    +g_{M}({\mathcal{A}}_{X_{j}}X_{j},{\mathcal{A}}_{X}X)\right).
\end{aligned}
$$
\begin{equation}
  Ric^*(X) =\frac{c}{4} \Bigg\{ (\ell -1) 
    +\sum_{j=1}^{\ell} \sum_{\alpha = 1}^{3} 
        3g_{M}(C_{J_{\alpha}}X, X_{j})^2
      \Bigg\}-  \sum_{j=1}^{\ell}3 \lVert{\mathcal{A}}_{X}X_{j} \rVert ^{2}.
      \label{eq-7.}
\end{equation}
\begin{equation*} 
Ric^*(X) \leq \frac{c}{4} \Bigg\{ (\ell -1) 
    +\sum_{j=1}^{\ell} \sum_{\alpha = 1}^{3} 
        3g_{M}(C_{J_{\alpha}}X, X_{j})^2
      \Bigg\}.\end{equation*}
By using (\ref{eq-7.}), we can say that the equality case of \rm(\ref{eq-7}) holds if and only if ${\mathcal{H}}(M)$ is integrable.
\begin{thm}
Let $\pi: M(c) \rightarrow N $ be an anti-invariant submersion from a quaternionic space form $(M(c),g_{M})$ onto a Riemannian manifold $(N,g_{N}).$ Then
\begin{equation}\label{eq-3.11}
\tau^* \leq \frac{c}{4}\Bigg\{ \ell (\ell -1) 
    +\sum_{i=1}^{\ell} \sum_{\alpha = 1}^{3} 
        3\lVert C_{J_{\alpha}}X_i\rVert^2
      \Bigg\}. 
\end{equation}
The equality case of \rm (\ref{eq-3.11}) holds if and only if ${\mathcal{H}}(M)$ is integrable.
\end{thm}
{\bf Proof:} 
We know that 
\[ {\tau}^{*} = \sum_{i=1}^{\ell} {Ric}^*(X_{i},X_{i}).\]
Using \rm(\ref{eq-7.}), we get
\begin{equation*}
\tau^* =\frac{c}{4}\sum_{i=1} ^{\ell}\Bigg\{ (\ell -1) 
    +\sum_{j=1}^{\ell} \sum_{\alpha = 1}^{3} 
        3g_{M}(C_{J_{\alpha}}X_i, X_{j})^2
      \Bigg\}-  \sum_{i,j=1}^{\ell}3 \lVert{\mathcal{A}}_{X_i}X_{j} \rVert ^{2}. 
\end{equation*}
\begin{equation}
\tau^* =\frac{c}{4}\Bigg\{\ell (\ell -1) 
    +\sum_{i=1}^{\ell} \sum_{\alpha = 1}^{3} 
        3\lVert C_{J_{\alpha}}X_i\rVert^2
      \Bigg\}-  \sum_{i,j=1}^{\ell}3 \lVert{\mathcal{A}}_{X_i}X_{j} \rVert ^{2}, 
      \label{eq-3.10}
\end{equation}
where $$\lVert C_{J_{\alpha}}X_i\rVert^2=\sum_{j=1}^{\ell}g_{M}(C_{J_{\alpha}}X_i, X_{j})^2.$$

\begin{equation*}
\tau^* \leq \frac{c}{4}\Bigg\{ \ell (\ell -1) 
    +\sum_{i=1}^{\ell} \sum_{\alpha = 1}^{3} 
        3\lVert C_{J_{\alpha}}X_i\rVert^2
      \Bigg\}. 
\end{equation*}
\vspace{1cm}

\noindent Define \rm{\cite{Gulbahar}} 
\begin{equation} \label{eq-3.12}
 {\mathcal{T}}_{ij}^{s}= g_{M}({\mathcal{T}}_{U_{i}}U_{j},X_{s}),
\end{equation}
\noindent where $ 1\leq i, j\leq r$ and $1 \leq s \leq \ell$
and
\begin{equation} \label{eq-3.13}
{\mathcal{A}}_{ij}^{\beta}=g_{M}({\mathcal{A}}_{X_{i}}X_{j},U_{\beta}),
\end{equation}
 where $1 \leq i, j \leq \ell $ and $1 \leq \beta \leq r.$  
\newline 
\noindent Since $\{U_1, \ldots, U_r\}$ is an orthonormal frame of ${\mathcal{V}}(M)$. Then, it follows that
 \[
 g(\nabla_E N, X) = \sum_{j=1}^r g((\nabla_E {\mathcal{T}})(U_j, U_j), X),
\]
 for any $E \in \Gamma (TM)$ and $X \in {\mathcal{H}}(M)$ \rm{\cite{Falcitelli-Ianus and Pastore}}.

The horizontal divergence of any vector field $X$ on ${\mathcal{H}}(M)$ is given by ${\delta}(X)$ and defined by
\[
{\delta}(X) = \sum_{i=1}^\ell g(\nabla_{X_i} X, X_i),
\]
where $\{X_1, \ldots, X_\ell \}$ is a local orthonormal frame of ${\mathcal{H}}(M)$. Then one has \rm{\cite{Besse}}
\begin{equation} \label{eq-3.14}
\delta(N)=\sum_{i=1}^{\ell} \sum_{k=1}^{r} g_{M}((\nabla _{X_{i}}{\mathcal{T}})(U_{k},U_{k}),X_{i}).
\end{equation}

\begin{thm}
Let $\pi: M(c) \rightarrow N $ be an anti-invariant submersion from a quaternionic space form $(M(c),g_{M})$ onto a Riemannian manifold $(N,g_{N}).$ Then we have 
\begin{equation} \label{eq-3.16}
\widehat{Ric}(U_{1}) \geq \frac{c}{4}(r-1) -\frac{1}{4}r^{2} \lVert H \rVert ^{2}.
\end{equation}
The equality case of {\rm (\ref{eq-3.16})} holds if and only if
$$
\begin{aligned}
{\mathcal{T}}_{11}^{s}&={\mathcal{T}}_{22}^{s}+ \cdots + {\mathcal{T}}_{rr}^{s},\\
{\mathcal{T}}_{ij}^{s}&=0, \quad j=2,\ldots, r.
\end{aligned}
$$
\end{thm}
{\bf Proof:} Using \rm(\ref{eq-3.12}) in \rm(\ref{eq-8.}), we have
\begin{equation} \label{eq-3.17}
 \hat{\tau}=\frac{c}{4}r(r-1)-r^{2}\lVert H \rVert ^{2} +\sum_{s=1}^{\ell} \sum _{i,j=1}^{r} ({\mathcal{T}}_{ij}^{s})^{2}.
\end{equation}
Thus using \rm(\ref{eq-3.15}) in \rm (\ref{eq-3.17}) can be written as
\begin{equation} \label{eq-3.18}
 \hat{\tau}=\frac{c}{4}r(r-1)- \frac{1}{2}r^{2}\lVert H \rVert ^{2}+({\mathcal{T}}_{11}^{s}-{\mathcal{T}}_{22}^{s}- ... - {\mathcal{T}}_{rr}^{s})^{2} + 2\sum_{s=1}^{\ell} \sum_{j=2}^{r} ({\mathcal{T}}_{1j}^{s})^{2}-2\sum_{s=1}^{\ell} \sum_{2\leq i<j \leq r} ({\mathcal{T}}_{ii}^{s}{\mathcal{T}}_{jj}^{s}-({\mathcal{T}}_{ij}^{s})^{2}).
\end{equation}
Then from \rm(\ref{eq-3.18}), we have
 \begin{equation} \label{eq-3.19}
 \hat{\tau} \geq\frac{c}{4}r(r-1)- \frac{1}{2}r^{2}\lVert H \rVert ^{2}-2\sum_{s=1}^{\ell} \sum_{2\leq i<j \leq r} ({\mathcal{T}}_{ii}^{s}{\mathcal{T}}_{jj}^{s}-({\mathcal{T}}_{ij}^{s})^{2}).
\end{equation}
Taking $U=W=U_{i}, \quad V=F=U_{j}$ in {\rm(\ref{eq-1})} and using {\rm(\ref{eq-3.12})}, we get
\begin{equation}\label{eq-3.20}
 \sum_{2 \leq i < j \leq r}R(U_{i},U_{j},U_{j},U_{i})=  \sum_{2 \leq i < j \leq r} \hat{R} (U_{i},U_{j},U_{j},U_{i})+ \sum_{s=1}^{\ell} \sum_{2\leq i<j \leq r} ({\mathcal{T}}_{ii}^{s}{\mathcal{T}}_{jj}^{s}-({\mathcal{T}}_{ij}^{s})^{2}).
\end{equation}
Using {\rm(\ref{eq-3.20})} in {\rm(\ref{eq-3.19})}, we get
\begin{equation} \label{eq-3.21}
 \hat{\tau}\geq\frac{c}{4}r(r-1)- \frac{1}{2}r^{2}\lVert H \rVert ^{2}+2 \sum_{2 \leq i < j \leq r} \hat{R} (U_{i},U_{j},U_{j},U_{i})-2 \sum_{2 \leq i < j \leq r}R(U_{i},U_{j},U_{j},U_{i}).
\end{equation}
We also have
\begin{equation} \label{eq-3.22}
 \hat{\tau}= 2\sum_{2 \leq i < j \leq r} \hat{R} (U_{i},U_{j},U_{j},U_{i})+ 2\sum_{j=1}^{r} \hat{R} (U_{1},U_{j},U_{j},U_{1}).
\end{equation}
Considering {\rm(\ref{eq-3.22})} in {\rm(\ref{eq-3.21})}, we obtain
\begin{equation}\label{eq-3.33}
 2\widehat{Ric}(U_{1}) \geq \frac{c}{4}r(r-1)- \frac{1}{2}r^{2}\lVert H \rVert ^{2}-2 \sum_{2 \leq i < j \leq r}R(U_{i},U_{j},U_{j},U_{i}).
\end{equation}
By {\rm(\ref{eq-2.9})}, we have
\begin{equation}\label{eq-24}
\sum_{2 \leq i < j \leq r}R(U_{i},U_{j},U_{j},U_{i})=\frac{c}{8}(r-1)(r-2).
\end{equation}
From {\rm(\ref{eq-3.33})} and {\rm(\ref{eq-24})}, we obtain  
\[\widehat{Ric}(U_{1}) \geq \frac{c}{4}(r-1)-\frac{1}{4}r^{2} \lVert H \rVert ^{2}.\]
\begin{thm}
    Let $\pi: M(c) \rightarrow N$ be an anti-invariant submersion from a quaternionic space form $ (M(c),g_{M})$ onto a Riemannian manifold $(N,g_{N}).$ Then we have 
    \begin{equation} \label{eq-3.25}
2 Ric^{*}(X_{1}) \leq \frac{c}{4}\left(2(\ell-1)
+3\lVert C_{J_\alpha}X_{1} \rVert ^{2} \right).
    \end{equation}
    The equality case of {\rm(\ref{eq-3.25})} holds if and only if 
    \[{\mathcal{A}}_{1j}^{\beta}=0, \hspace{0.1mm} j=2,...,\ell.\]
\end{thm}
{\bf Proof:}  Using {\rm(\ref{eq-3.13})} in {\rm(\ref{eq-3.10})}, we can write
\begin{equation} \label{eq-3.26}
\tau_{*}= \frac{c}{4}\left( \ell(\ell-1) + 3 \sum_{i=1}^{\ell}\sum_{\alpha=1}^{3}\lVert C_{J_\alpha }X_i \rVert ^{2}\right)-3 \sum_{\beta=1}^{r}\sum_{i,j=1}^{\ell}({\mathcal{A}}_{ij}^{\beta})^{2}.
\end{equation}
Using that ${\mathcal{A}}$ is anti-symmetric on ${\mathcal{H}}(M(c)),$ {\rm(\ref{eq-3.26})} can be written as 
\begin{equation} \label{eq-3.27}
\tau_{*}= \frac{c}{4}\left( \ell(\ell-1) + 3 \sum_{i=1}^{\ell}\sum_{\alpha=1}^{3}\lVert C_{J_\alpha }X_i \rVert ^{2}\right)-6\sum_{\beta=1}^{r}\sum_{j=2}^{\ell}({\mathcal{A}}_{1j}^{\beta})^{2}-6\sum_{\beta=1}^{r}\sum_{2 \leq i<j \leq \ell}({\mathcal{A}}_{ij}^{\beta})^{2}.
\end{equation}
Taking $X=H=X_{i},\quad Y=Z=X_{j}$ in {\rm(\ref{eq-2})} and using {\rm(\ref{eq-3.13})}, we obtain
\begin{equation} \label{eq-3.28}
 2\sum_{2 \leq i<j \leq \ell}R(X_{i},X_{j},X_{j},X_{i})=2\sum_{2 \leq i<j \leq \ell}R^{*}(X_{i},X_{j},X_{j},X_{i})+6\sum_{\beta=1}^{r}\sum_{2 \leq i<j \leq \ell}({\mathcal{A}}_{ij}^{\beta})^{2}.
\end{equation}
Using {\rm(\ref{eq-3.28})} in {\rm(\ref{eq-3.27})}, we get
\begin{equation} \label{eq-3.29}
\begin{aligned}
\tau_{*}&= \frac{c}{4}\left( \ell(\ell-1) + 3 \sum_{i=1}^{\ell}\sum_{\alpha=1}^{3}\lVert C_{J_\alpha }X_i \rVert ^{2}\right)-6\sum_{\beta=1}^{r}\sum_{j=2}^{\ell}({\mathcal{A}}_{1j}^{\beta})^{2}\\
&+2\sum_{2 \leq i<j \leq \ell}R^{*}(X_{i},X_{j},X_{j},X_{i})-2\sum_{2 \leq i<j \leq \ell}R(X_{i},X_{j},X_{j},X_{i}).
\end{aligned}
\end{equation}
From {\rm(\ref{eq-2.9})}, we obtain
\begin{equation} \label{eq-3.30}
\sum_{2 \leq i<j \leq \ell}R(X_{i},X_{j},X_{j},X_{i})=\frac{c}{4}\left( \frac{(\ell-1)(\ell-2)}{2}+ 3 \sum_{2 \leq i<j \leq \ell}\  \sum_{\alpha=1}^{3} g_{M}(C_{J_\alpha}X_{i},X_{j}) ^{2}\right).
\end{equation}
Then from {\rm(\ref{eq-3.29})} and {\rm(\ref{eq-3.30})}
\begin{equation} \label{eq-3.31}
 2Ric^{*}(X_{1})= \frac{c}{4}\left(2(\ell-1)
+3\lVert C_{J_\alpha}X_{1} \rVert ^{2} \right)- 6\sum_{\beta=1}^{r}\sum_{j=2}^{\ell}({\mathcal{A}}_{1j}^{\beta})^{2}.
\end{equation}
Therefore, we get the inequality
\begin{equation*} 
2 Ric^{*}(X_{1}) \leq \frac{c}{4}\left(2(\ell-1)
+3\lVert C_{J_\alpha}X_{1} \rVert ^{2} \right).
\end{equation*}
Now, we compute the Chen-Ricci inequality between vertical and horizontal distributions. For the scalar curvature $\tau$ of $M(c)$, we can write
\begin{equation}
\begin{aligned}
 \tau & =\sum_{s=1}^\ell \operatorname{Ric}\left(X_s, X_s\right)+\sum_{k=1}^r \operatorname{Ric}\left(U_k, U_k\right), 
\end{aligned}
\end{equation}\begin{equation}
\begin{aligned}
\tau & =\sum_{j, k=1}^r R\left(U_j, U_k, U_k, U_j\right)+\sum_{i=1}^\ell \sum_{k=1}^r R\left(X_i, U_k, U_k, X_i\right) \\
& +\sum_{i, s=1}^\ell R\left(X_i, X_s, X_s, X_i\right)+\sum_{s=1}^\ell \sum_{j=1}^r R\left(U_j, X_s, X_s, U_j\right) .
\end{aligned}
\label{eq-9}
\end{equation}
Let us denote,
\begin{equation}\label{eq-35}
\begin{aligned}
& \left\|\mathcal{T}^{\mathcal{V}}\right\|^2=\sum_{i=1}^\ell \sum_{k=1}^r g_M\left({\mathcal{T}}_{U_k} X_i, {\mathcal{T}}_{U_k} X_i\right).
\end{aligned}
\end{equation}
\begin{equation}\label{eq-36}
\begin{aligned}
& \left\|\mathcal{T}^{\mathcal{H}}\right\|^2=\sum_{j, k=1}^r g_M\left({\mathcal{T}}_{U_j} U_{k,} {\mathcal{T}}_{U_j} U_k\right). 
\end{aligned}
\end{equation}
\begin{equation}\label{eq-37}
\begin{aligned}
& \left\|\mathcal{A}^{\mathcal{V}}\right\|^2=\sum_{i, j=1}^\ell g_M\left({\mathcal{A}}_{X_i} X_{j,} {\mathcal{A}}_{X_i} X_{j}\right).
\end{aligned}
\end{equation}
\begin{equation}\label{eq-38}
\begin{aligned}
& \left\|\mathcal{A}^{\mathcal{H}}\right\|^2=\sum_{i=1}^\ell \sum_{k=1}^r g_M\left({\mathcal{A}}_{X_i} U_{k,} {\mathcal {A}}_{X_i} U_{k}\right).
\end{aligned}
\end{equation}
\begin{thm}
Let $\pi: M(c) \rightarrow N$ be an anti-invariant submersion from a quaternionic space form $(M(c),g_{M})$ onto a Riemannian manifolds $(N,g_{N}).$ Then we have
\begin{equation} \label{eq-4}
\begin{aligned}
&\frac{c}{4}\left(\ell r+\ell+r+\sum_{\alpha =1}^{3}\left( \sum_{i=1}^{\ell} \lVert B_{J_\alpha}X_i \rVert ^{2}+3 \lVert C_{J_\alpha}X_{1} \rVert ^{2}\right) \right) \\
&\leq \widehat{Ric}(U_{1})+Ric^{*}(X_{1}) + \frac{1}{4} r^{2} \lVert H \rVert ^{2}+3 \sum_{\beta =1}^{r} \sum_{s=2}^{\ell} ({\mathcal{A}}_{1s}^{\beta})^{2} -\delta (N) + \left\| {\mathcal{T}}^{\mathcal{V}}\right\|^{2}-\left\| {\mathcal{A}}^{\mathcal{H}} \right\|^{2}.
\end{aligned}
\end{equation}
The equality case of {\rm(\ref{eq-4})} holds if and only if
$$
\begin{aligned}
{\mathcal{T}}_{11}^{s}&={\mathcal{T}}_{22}^{s}+ \cdots + {\mathcal{T}}_{rr}^{s},\\
{\mathcal{T}}_{ij}^{s}&=0,  \quad j=2, \cdots, r.
\end{aligned}
$$
\end{thm}
{\bf Proof:} Since $M(c)$ is a quaternionic space form, using {\rm(\ref{eq-9})} we obtain
\begin{equation} \label{eq-3.34}
 \tau= \frac{c}{4} \left[ (\ell+r)(\ell+r-1) +3 \sum_{\alpha=1}^{3}\sum_{i=1}^{\ell} \left(  \lVert C_{J_\alpha}X_i \rVert ^{2} + 2\sum_{k=1}^{r}g_{M}(B_{J_\alpha}X_{i},U_{k})^{2}\right)\right].
\end{equation}
Now, we define 
\begin{equation} \label{eq-3.35}
\lVert B_{J_\alpha}X_{i} \rVert ^{2}= \sum_{k=1}^{r} g_{M}(B_{J_\alpha}X_{i},U_{k})^{2}.
\end{equation}
On the other hand, using Gauss-Codazzi type equations {\rm(\ref{eq-1})}, {\rm(\ref{eq-2})} and {\rm(\ref{eq-3})}, we obtain
\begin{equation}\label{eq-3.36}
\begin{aligned}
 \tau &=   \hat{\tau} + \tau^{*}+r^{2} \lVert H \rVert ^{2}- \sum_{k,j=1}^{r}g_{M}({\mathcal{T}}_{U_{k}}U_{j},{\mathcal{T}}_{U_{k}}U_{j})+3\sum_{i,s=1}^{\ell}g_{M}({\mathcal{A}}_{X_{i}}X_{s}, {\mathcal{A}}_{X_{i}}X_{s})\\
& -\sum_{i=1}^{\ell}\sum_{k=1}^{r}g_{M}\left( (\nabla_{X_{i}}{\mathcal{T}})(U_{k},U_{k}),X_{i}\right) + \sum_{i=1}^{\ell}\sum_{k=1}^{r}\left( g_{M}({\mathcal{T}}_{U_{k}}X_{i},{\mathcal{T}}_{U_{k}}X_{i})-g_{M}({\mathcal{A}}_{X_{i}}U_{k}, {\mathcal{A}}_{X_{i}}U_{k}) \right)\\
&-\sum_{s=1}^{\ell}\sum_{j=1}^{r}g_{M}\left( (\nabla_{X_{s}}{\mathcal{T}})(U_{i},U_{j}),X_{s}\right)+\sum_{s=1}^{\ell}\sum_{j=1}^{r}\left( g_{M}({\mathcal{T}}_{U_{j}}X_{s},{\mathcal{T}}_{U_{j}}X_{s})-g_{M}({\mathcal{A}}_{X_{s}}U_{j}, {\mathcal{A}}_{X_{s}}U_{j}) \right).
\end{aligned}
\end{equation}
Then using {\rm(\ref{eq-3.13})}, {\rm(\ref{eq-3.14})}, {\rm(\ref{eq-3.15})} and {\rm(\ref{eq-3.36})}, we obtain
\begin{equation}\label{eq-3.37}
\begin{aligned}
 \tau &= \hat{\tau} + \tau^{*}+\frac{1}{2}r^{2} \lVert H \rVert ^{2}-({\mathcal{T}}_{11}^{s}-{\mathcal{T}}_{22}^{s}- \cdots- {\mathcal{T}}_{rr}^{s})^{2}-2\sum_{s=1}^{\ell}\sum_{j=2}^{r}({\mathcal{T}}_{1j}^{s})^{2}\\
&+2\sum_{s=1}^{\ell}\sum_{2 \leq i<j \leq r}({\mathcal{T}}_{ii}^{s}{\mathcal{T}}_{jj}^{s}-({\mathcal{T}}_{ij}^{s})^{2})+6\sum_{\beta=1}^{r}\sum_{s=2}^{\ell}({\mathcal{A}}_{1s}^{\beta})^{2}+6\sum_{\beta=1}^{r}\sum_{2 \leq i<j \leq \ell}({\mathcal{A}}_{is}^{\beta})^{2}\\
&+ \sum_{i=1}^{\ell}\sum_{k=1}^{r}\left( g_{M}({\mathcal{T}}_{U_{k}}X_{i},{\mathcal{T}}_{U_{k}}X_{i})-g_{M}({\mathcal{A}}_{X_{i}}U_{k}, {\mathcal{A}}_{X_{i}}U_{k}) \right)-2\delta(N)\\
&+\sum_{s=1}^{\ell}\sum_{j=1}^{r}\left( g_{M}({\mathcal{T}}_{U_{j}}X_{s},{\mathcal{T}}_{U_{j}}X_{s})-g_{M}({\mathcal{A}}_{X_{s}}U_{j}, {\mathcal{A}}_{X_{s}}U_{j}) \right).
\end{aligned}
\end{equation}
Using {\rm(\ref{eq-3.20}),\rm(\ref{eq-3.28}),\rm(\ref{eq-3.34})} and {\rm(\ref{eq-3.35})} in {\rm(\ref{eq-3.37})}
\begin{eqnarray}\label{eq-3.38}
\frac{c}{4} \left[ (\ell+r)(\ell+r-1) + 3 \sum_{\alpha=1}^{3}\sum_{i=1}^{\ell} \left(  \lVert C_{J_\alpha}X_i \rVert ^{2} + 2\lVert B_{J_\alpha}X_i \rVert ^{2}\right) \right]
= \hat{\tau} + \tau^{*}+\frac{1}{2}r^{2} \lVert H \rVert ^{2}\notag\\
-({\mathcal{T}}_{11}^{s}-{\mathcal{T}}_{22}^{s}- \cdots- {\mathcal{T}}_{rr}^{s})^{2}
+2\sum_{2 \leq i < j \leq r}R(U_{i},U_{j},U_{j},U_{i})-2\sum_{2 \leq i < j \leq r}\hat{R}(U_{i},U_{j},U_{j},U_{i})\notag\\
-2\sum_{s=1}^{\ell}\sum_{j=2}^{r}({\mathcal{T}}_{1j}^{s})^{2}
+6\sum_{\beta=1}^{r}\sum_{s=2}^{\ell}({\mathcal{A}}_{1s}^{\beta})^{2}
+2\sum_{2 \leq i<j \leq \ell}R(X_{i},X_{j},X_{j},X_{i})-2\sum_{2 \leq i<j \leq \ell}R^{*}(X_{i},X_{j},X_{j},X_{i})\notag\\ + 2\sum_{i=1}^{\ell}\sum_{k=1}^{r}( g_{M}({\mathcal{T}}_{U_{k}}X_{i},{\mathcal{T}}_{U_{k}}X_{i})
-g_{M}({\mathcal{A}}_{X_{i}}U_{k}, {\mathcal{A}}_{X_{i}}U_{k}) )-2\delta(N).\quad
\end{eqnarray}
\begin{equation}\label{eq-3.39}
\begin{aligned}
&\frac{c}{4} \left[ (\ell+r)(\ell+r-1) + 3 \sum_{\alpha=1}^{3}\sum_{i=1}^{\ell} \left(  \lVert C_{J_\alpha}X_i \rVert ^{2} + 2\lVert B_{J_\alpha}X_i \rVert ^{2}\right) \right]\\&=2\widehat{Ric}(U_{1}) + 2{Ric}^{*}(X_{1})+2\sum_{2 \leq i<j \leq \ell}R(X_{i},X_{j},X_{j},X_{i}) \\
&+\frac{1}{2}r^{2} \lVert H \rVert ^{2}-({\mathcal{T}}_{11}^{s}-{\mathcal{T}}_{22}^{s}- \cdots- {\mathcal{T}}_{rr}^{s})^{2}
-2\sum_{s=1}^{\ell}\sum_{j=2}^{r}({\mathcal{T}}_{1j}^{s})^{2}+6\sum_{\beta=1}^{r}\sum_{s=2}^{\ell}({\mathcal{A}}_{1s}^{\beta})^{2}-2\delta(N)\\
&+ 2\sum_{i=1}^{\ell}\sum_{k=1}^{r}\left( g_{M}({\mathcal{T}}_{U_{k}}X_{i},{\mathcal{T}}_{U_{k}}X_{i})
-g_{M}({\mathcal{A}}_{X_{i}}U_{k}, {\mathcal{A}}_{X_{i}}U_{k}) \right)+2\sum_{2 \leq i < j \leq r}R(U_{i},U_{j},U_{j},U_{i}).
\end{aligned}
\end{equation}
\begin{equation}\label{eq-3.39.}
\begin{aligned}
&\frac{c}{4} \left[ (\ell+r)(\ell+r-1) + 3 \sum_{\alpha=1}^{3}\sum_{i=1}^{\ell} \left(  \lVert C_{J_\alpha}X_i \rVert ^{2} + 2\lVert B_{J_\alpha}X_i \rVert ^{2}\right) \right]\\&=2\widehat{Ric}(U_{1}) + 2{Ric}^{*}(X_{1})+\frac{c}{2} \left( \frac{(\ell-1)(\ell-2)}{2}+ 3 \sum_{\alpha=1}^{3}\sum_{2 \leq i<j \leq \ell}^{} g_M( C_{J_\alpha}X_i, X_j) ^{2} \right) \\
&+\frac{1}{2}r^{2} \lVert H \rVert ^{2}-({\mathcal{T}}_{11}^{s}-{\mathcal{T}}_{22}^{s}- \cdots- {\mathcal{T}}_{rr}^{s})^{2}
-2\sum_{s=1}^{\ell}\sum_{j=2}^{r}({\mathcal{T}}_{1j}^{s})^{2}+6\sum_{\beta=1}^{r}\sum_{s=2}^{\ell}({\mathcal{A}}_{1s}^{\beta})^{2}-2\delta(N)\\
&+ 2\sum_{i=1}^{\ell}\sum_{k=1}^{r}\left( g_{M}({\mathcal{T}}_{U_{k}}X_{i},{\mathcal{T}}_{U_{k}}X_{i})
-g_{M}({\mathcal{A}}_{X_{i}}U_{k}, {\mathcal{A}}_{X_{i}}U_{k}) \right)+\frac{c}{4}(r-1)(r-2).
\end{aligned}
\end{equation}
Using {\rm(\ref{eq-24}), \rm(\ref{eq-3.30}), \rm(\ref{eq-35})} and {\rm(\ref{eq-38})} in {\rm(\ref{eq-3.39})}, we obtain the inequality
$$
\begin{aligned}
&\frac{c}{4}\left(\ell r+\ell+r+\sum_{\alpha =1}^{3}\left( \sum_{i=1}^{\ell} \lVert B_{J_\alpha}X_i \rVert ^{2}+3 \lVert C_{J_\alpha}X_{1} \rVert ^{2}\right) \right) \\
&\leq \widehat{Ric}(U_{1})+Ric^{*}(X_{1}) + \frac{1}{4} r^{2} \lVert H \rVert ^{2}+3 \sum_{\beta =1}^{r} \sum_{s=2}^{\ell} ({\mathcal{A}}_{1s}^{\beta})^{2} -\delta (N) + \left\| {\mathcal{T}}^{\mathcal{V}}\right\|^{2}-\left\| {\mathcal{A}}^{\mathcal{H}} \right\|^{2}.
\end{aligned}
$$

\vspace{0.3cm}
\noindent From {\rm(\ref{eq-35}), \rm(\ref{eq-38}), \rm(\ref{eq-3.34}), \rm(\ref{eq-3.35})} and {\rm(\ref{eq-3.36})}, we obtain
\begin{equation}\label{eq-3.47}
\begin{aligned}
& \frac{c}{4}\left((\ell+r)(\ell+r-1)+\sum_{\alpha =1}^{3}\left( \sum_{i=1}^{\ell} 6\lVert B_{J_\alpha}X_i \rVert ^{2}+3 \lVert C_{J_\alpha}X_{1} \rVert ^{2}\right) \right)\\&= \hat{\tau}+ \tau^* 
 +r^2\|H\|^2-\left\|{\mathcal{T}}^{\mathcal{H}}\right\|^2+3\left\|{\mathcal{A}}^{\mathcal{V}}\right\|^2-2 \delta(N)+2\left\|{\mathcal{T}}^{\mathcal{V}}\right\|^2-2\left\|{\mathcal{A}}^{\mathcal{H}}\right\|^2.
\end{aligned}
\end{equation}
From {\rm(\ref{eq-3.47})}, we obtain the following theorem:
\begin{thm} \label{th-7}
Let $\pi: M(c) \rightarrow N$ be an anti-invariant submersion from a quaternionic space form $\left(M(c), g_M\right)$ onto a Riemannian manifold $(N,g_{N})$. Then we have
\begin{equation} \label{eq-3.48}
\begin{aligned}
 \hat{\tau}+ \tau^* & \leq  \frac{c}{4}\left((\ell+r)(\ell+r-1)+\sum_{\alpha =1}^{3}\left( \sum_{i=1}^{\ell} 6\lVert B_{J_\alpha}X_i \rVert ^{2}+3 \lVert C_{J_\alpha}X_{1} \rVert ^{2}\right) \right) \\
&-r^2\|H\|^2+\left\|{\mathcal{T}}^{\mathcal{H}}\right\|^2+2 \delta(N)-2\left\|{\mathcal{T}}^{\mathcal{V}}\right\|^2+2\left\|{\mathcal{A}}^{\mathcal{H}}\right\|^2,
\end{aligned}
\end{equation}
\begin{equation}\label{eq-3.49}
\begin{aligned}
 \hat{\tau}+ \tau^* & \geq \frac{c}{4}\left((\ell+r)(\ell+r-1)+\sum_{\alpha =1}^{3}\left( \sum_{i=1}^{\ell} 6\lVert B_{J_\alpha}X_i \rVert ^{2}+3 \lVert C_{J_\alpha}X_{1} \rVert ^{2}\right) \right) \\
& -r^2\|H\|^2+\left\|{\mathcal{T}}^{\mathcal{H}}\right\|^2-3\left\|{\mathcal{A}}^{\mathcal{V}}\right\|^2+2 \delta(N)-2\left\|{\mathcal{T}}^{\mathcal{V}}\right\|^2.
\end{aligned}
\end{equation}
\noindent Equality cases of {\rm(\ref{eq-3.48})} and {\rm(\ref{eq-3.49})}  hold for all $p \in M$ if and only if horizontal distribution $\mathcal{H}$ is integrable.
\end{thm}
From Theorem~\ref{th-7}, we have the following corollary:
\begin{cor}Let $\pi: M(c) \rightarrow N$ be an anti-invariant submersion from a quaternionic space form  $(M(c), g_M)$  onto a Riemannian manifold  $(N, g_N)$  such that each fibres be totally geodesic. Then we have
\begin{equation}\label{eq-3.50}
\begin{aligned}
&  \hat{\tau}+ \tau^* \leq \frac{c}{4}\left((\ell+r)(\ell+r-1)+\sum_{\alpha =1}^{3}\left( \sum_{i=1}^{\ell} 6\lVert B_{J_\alpha}X_i \rVert ^{2}+3 \lVert C_{J_\alpha}X_{1} \rVert ^{2}\right) \right)+2\left\|{\mathcal{A}}^{\mathcal{H}}\right\|^2, 
\end{aligned}
\end{equation}
\begin{equation}\label{eq-3.51}
\begin{aligned}
&  \hat{\tau}+ \tau^* \geq \frac{c}{4}\left((\ell+r)(\ell+r-1)+\sum_{\alpha =1}^{3}\left( \sum_{i=1}^{\ell} 6\lVert B_{J_\alpha}X_i \rVert ^{2}+3 \lVert C_{J_\alpha}X_{1} \rVert ^{2}\right) \right)-3\left\|{\mathcal{A}}^{\mathcal{V}}\right\|^2.
\end{aligned}
\end{equation}
\noindent Equality cases of {\rm(\ref{eq-3.50})} and {\rm(\ref{eq-3.51})} hold for all $p \in M$ if and only if horizontal distribution $\mathcal{H}$ is integrable.
\end{cor}
From {\rm(\ref{eq-3.47})} we obtain following theorem:
\begin{thm} \label{th-9} Let $\pi: M(c) \rightarrow N$ be an anti-invariant submersion from a quaternionic space form  $(M(c), g_M)$  onto a Riemannian manifold  $(N,g_{N}).$ Then we have
\begin{equation}\label{eq-3.52}
\begin{aligned}
 \hat{\tau}+ \tau^* & \geq \frac{c}{4}\left((\ell+r)(\ell+r-1)+\sum_{\alpha =1}^{3}\left( \sum_{i=1}^{\ell} 6\lVert B_{J_\alpha}X_i \rVert ^{2}+3 \lVert C_{J_\alpha}X_{1} \rVert ^{2}\right) \right) \\
& -r^2\|H\|^2+2 \delta(N)-2\left\|{\mathcal{T}}^{\mathcal{V}}\right\|^2+2\left\|{\mathcal{A}}^{\mathcal{H}}\right\|^2-3\left\|{\mathcal{A}}^{\mathcal{V}}\right\|^2,
\end{aligned}
\end{equation}
\begin{equation}\label{eq-3.53}
\begin{aligned}
 \hat{\tau}+ \tau^* & \leq \frac{c}{4}\left((\ell+r)(\ell+r-1)+\sum_{\alpha =1}^{3}\left( \sum_{i=1}^{\ell} 6\lVert B_{J_\alpha}X_i \rVert ^{2}+3 \lVert C_{J_\alpha}X_{1} \rVert ^{2}\right) \right) \\
& -r^2\|H\|^2+\left\|{\mathcal{T}}^{\mathcal{H}}\right\|^2+2 \delta(N)+2\left\|{\mathcal{A}}^{\mathcal{H}}\right\|^2-3\left\|{\mathcal{A}}^{\mathcal{V}}\right\|^2.
\end{aligned}
\end{equation}
\noindent Equality cases of {\rm(\ref{eq-3.52})} and {\rm(\ref{eq-3.53})} hold for all $p \in M$ if and only if the fibre through $p$ of $\pi$ is a totally geodesic submanifold of $M$.
\end{thm}
From Theorem~\ref{th-9}, we have the following corollary:
\begin{cor}
Let $\pi:$ Let $\pi: M(c) \rightarrow N$ be an anti-invariant submersion from a quaternionic space form  $(M(c), g_M)$  onto a Riemannian manifold  $(N, g_N)$  such that ${\mathcal{H}}$ be integrable. Then we have
\begin{equation}\label{eq-3.54}
\begin{aligned}
 \hat{\tau}+ \tau^* & \geq \frac{c}{4}\left((\ell+r)(\ell+r-1)+\sum_{\alpha =1}^{3}\left( \sum_{i=1}^{\ell} 6\lVert B_{J_\alpha}X_i \rVert ^{2}+3 \lVert C_{J_\alpha}X_{1} \rVert ^{2}\right) \right) \\ & -r^2\|H\|^2+2 \delta(N)-2\left\|{\mathcal{T}}^{\mathcal{V}}\right\|^2,
\end{aligned}
\end{equation}
\begin{equation}\label{eq-3.55}
\begin{aligned}
 \hat{\tau}+ \tau^* & \leq \frac{c}{4}\left((\ell+r)(\ell+r-1)+\sum_{\alpha =1}^{3}\left( \sum_{i=1}^{\ell} 6\lVert B_{J_\alpha}X_i \rVert ^{2}+3 \lVert C_{J_\alpha}X_{1} \rVert ^{2}\right) \right) \\
& -r^2\|H\|^2+2 \delta(N)+\left\|{\mathcal{T}}^{\mathcal{H}}\right\|^2.
\end{aligned}
\end{equation}
\noindent Equality cases of {\rm(\ref{eq-3.54})} and {\rm(\ref{eq-3.55})} hold for all $p \in M$ if and only if the fibre through $p$ of $\pi$ is a totally geodesic submanifold of $M$.
\end{cor}
Using Lemma~\ref{lem-11} in {\rm(\ref{eq-3.47})}, we obtain the following results:
\begin{thm} Let $\pi: M(c) \rightarrow N$ be an anti-invariant submersion from a quaternionic space form $\left(M(c), g_M\right)$ onto a Riemannian manifold $(N,g_{N}).$ Then we have
\begin{equation}\label{eq-3.56}
\begin{aligned}
& \frac{c}{4}\left((\ell+r)(\ell+r-1)+\sum_{\alpha =1}^{3}\left( \sum_{i=1}^{\ell} 6\lVert B_{J_\alpha}X_i \rVert ^{2}+3 \lVert C_{J_\alpha}X_{1} \rVert ^{2}\right) \right)\\ & \leq \hat{\tau}+ \tau^* 
 +r^2\|H\|^2+2\left\|{\mathcal{T}}^{\mathcal{V}}\right\|^2+3\left\|{\mathcal{A}}^{\mathcal{V}}\right\|^2-2 \delta(N)-2 \sqrt{2}\left\|{\mathcal{A}} ^{\mathcal{H}}\right\|\left\|{\mathcal{T}}^{\mathcal{H}}\right\|.
\end{aligned}
\end{equation}
\noindent Equality cases of {\rm(\ref{eq-3.56})} hold for all $p \in M$ if and only if $\left\|{\mathcal{A}}^{\mathcal{H}}\right\|=\left\|{\mathcal{T}}^{\mathcal{H}}\right\|$.
\end{thm}
\begin{thm} Let $\pi: M(c) \rightarrow N$ be an anti-invariant submersion from a quaternionic space form $\left(M(c), g_M\right)$ onto a Riemannian manifold  $(N, g_N).$  Then we have
\begin{equation}\label{eq-3.57}
\begin{aligned}
& \frac{c}{4}\left((\ell+r)(\ell+r-1)+\sum_{\alpha =1}^{3}\left( \sum_{i=1}^{\ell} 6\lVert B_{J_\alpha}X_i \rVert ^{2}+3 \lVert C_{J_\alpha}X_{1} \rVert ^{2}\right) \right) \geq \\
&  \hat{\tau}+ \tau^*+r^2\|H\|^2-\left\|{\mathcal{T}}^{\mathcal{H}}\right\|^2-2 \delta(N)-2\left\|{\mathcal{A}}^{\mathcal{H}}\right\|^2+2 \sqrt{6}\left\|{\mathcal{A}}^{\mathcal{V}}\right\|\left\|{\mathcal{T}}^{\mathcal{V}}\right\|.
\end{aligned}
\end{equation}
Equality cases of {\rm(\ref{eq-3.57})} hold for all $p \in M$ if and only if $\left\|{\mathcal{A}}^{\mathcal{V}}\right\|=\left\|{\mathcal{T}}^{\mathcal{V}}\right\|$.
\end{thm}
\begin{thm} \label{th-15}Let $\pi: M(c) \rightarrow N$ be an anti-invariant submersion from a quaternionic space form $\left(M(c), g_M\right)$ onto a Riemannian manifold  $(N, g_N).$ Then we have
\begin{equation}\label{eq-3.58}
\begin{aligned}
& \frac{c}{4}\left((\ell+r)(\ell+r-1)+\sum_{\alpha =1}^{3}\left( \sum_{i=1}^{\ell} 6\lVert B_{J_\alpha}X_i \rVert ^{2}+3 \lVert C_{J_\alpha}X_{1} \rVert ^{2}\right) \right) \\ & \leq  \hat{\tau}+ \tau^* 
+r(r-1)\|H\|^2+3\left\|{\mathcal{A}}^{\mathcal{V}}\right\|^2-2 \delta(N)+2\left\|{\mathcal{T}}^{\mathcal{V}}\right\|^2-2\left\|{\mathcal{A}}^{\mathcal{H}}\right\|^2.
\end{aligned}
\end{equation}
\noindent Equality case of {\rm(\ref{eq-3.58})} holds for all $p \in M$ if and only if we have the following statements:
\begin{enumerate}
    \item $\pi$ is a Riemannian submersion that has a totally umbilical fibres.
    \item ${\mathcal{T}}_{i j}=0$, for $i \neq j \in\{1,2, \ldots, r\}$.
\end{enumerate} 
\end{thm}
{\bf Proof:} From {\rm(\ref{eq-3.47})}, we have
\begin{equation} \label{eq-3.59}
\begin{aligned}
& \frac{c}{4}\left((\ell+r)(\ell+r-1)+\sum_{\alpha =1}^{3}\left( \sum_{i=1}^{\ell} 6\lVert B_{J_\alpha}X_i \rVert ^{2}+3 \lVert C_{J_\alpha}X_{1} \rVert ^{2}\right) \right) 
=\hat{\tau}+ \tau^*+r^2\|H\|^2 \\
&-\sum_{s=1}^\ell \sum_{j=1}^r\left({\mathcal{T}}_{j j}^s\right)^2
-2 \sum_{s=1}^\ell \sum_{j \neq k}^r\left({\mathcal{T}}_jk^s\right)^2+3\left\|{\mathcal{A}}^{\mathcal{V}}\right\|^2-2 \delta(N)+2\left\|{\mathcal{T}}^{\mathcal{V}}\right\|^2-2\left\|{\mathcal{A}}^{\mathcal{H}}\right\|^2 .
\end{aligned}
\end{equation}
Using Lemma~(\ref{lem-14}) in {\rm(\ref{eq-3.59})}, we obtain
$$
\begin{aligned}
&\frac{c}{4}\left((\ell+r)(\ell+r-1)+\sum_{\alpha =1}^{3}\left( \sum_{i=1}^{\ell} 6\lVert B_{J_\alpha}X_i \rVert ^{2}+3 \lVert C_{J_\alpha}X_{1} \rVert ^{2}\right) \right) 
\leq  \hat{\tau}+ \tau^*+r^2\|H\|^2\\
&-\frac{1}{r} \sum_{s=1}^\ell\left(\sum_{j=1}^r {\mathcal{T}}_{j j}^s\right)^2-2 \sum_{s=1}^\ell \sum_{j \neq k}^r\left({\mathcal{T}}_{jk}^s\right)^2+3\left\|{\mathcal { A }} ^ {\mathcal{V}}\right\|^2-2 \delta(N)+2\|{\mathcal{T}}^{\mathcal{V}}\|^2-2\left\|{\mathcal{\mathcal { A }} ^ { H }}\right\|^2,
\end{aligned}
$$
which is equivalent to {\rm (\ref{eq-3.58})}. Equality case of {\rm (\ref{eq-3.58})} holds for all $p \in M$ if and only if
$$
{\mathcal{T}}_{11}={\mathcal{T}}_{22}=\ldots={\mathcal{T}}_{rr} \quad \text { and } \quad \sum_{s=1}^\ell \sum_{j \neq k}^r\left({\mathcal{T}}_{jk}^s\right)^2=0,
$$
which completes proof of the theorem.
\newline \noindent
Using a similar proof approach as in Theorem~\ref{th-15}, we have the following theorem:
\begin{thm} \label{th-16} Let $\pi: M(c) \rightarrow N$ be an anti-invariant submersion from a quaternionic space form ( $M(c), g_M$ ) onto a Riemannian manifold $\left(N, g_N\right)$.Then we have
\begin{equation} \label{eq-3.60}
\begin{aligned}
& \frac{c}{4}\left((\ell+r)(\ell+r-1)+\sum_{\alpha =1}^{3}\left( \sum_{i=1}^{\ell} 6\lVert B_{J_\alpha}X_i \rVert ^{2}+3 \lVert C_{J_\alpha}X_{1} \rVert ^{2}\right) \right)  \geq  \hat{\tau}+ \tau^* \\
& +r^2\|H\|^2-\left\|{\mathcal{T}}^{\mathcal{H}}\right\|^2+\frac{3}{\ell}\operatorname{tr}\left({\mathcal{A}}^{\mathcal{V}}\right)^2-2 \delta(N)+2\left\|{\mathcal{T}}^{\mathcal{V}}\right\|^2-2\left\|{\mathcal{A}}^{\mathcal{H}}\right\|^2.
\end{aligned}
\end{equation}
Equality case of {\rm(\ref{eq-3.60})} holds for all $p \in M$ if and only if ${\mathcal{A}}_{11}={\mathcal{A}}_{22}=\ldots={\mathcal{A}}_{\ell \ell}$ and ${\mathcal{A}}_{i j}=0$, for $i \neq j \in\{1,2, \ldots, \ell \}$.
\end{thm}
From Theorem~\ref{th-16}, we have the following corollary:
\begin{cor}Let $\pi: M(c) \rightarrow N$ be an anti-invariant submersion from a quaternionic space form ( $M(c), g_M$ ) onto a Riemannian manifold $\left(N, g_N\right)$ such that each fiber is totally geodesic. Then we have
\begin{equation} \label{eq-3.61}
\begin{aligned}
 \frac{c}{4}\left((\ell+r)(\ell+r-1)+\sum_{\alpha =1}^{3}\left( \sum_{i=1}^{\ell} 6\lVert B_{J_\alpha}X_i \rVert ^{2}+3 \lVert C_{J_\alpha}X_{1} \rVert ^{2}\right) \right)  
 \geq \hat{\tau}+ \tau^*+\frac{3}{\ell} \operatorname{tr}\left({\mathcal{A}}^{\mathcal{V}}\right)^2-2\left\|{\mathcal{A}}^{\mathcal{H}}\right\|^2.
\end{aligned}
\end{equation}
Equality case of {\rm(\ref{eq-3.61})} holds for all $p \in M$ if and only if ${\mathcal{A}}_{11}={\mathcal{A}}_{22}=\ldots={\mathcal{A}}_{\ell \ell}$ and ${\mathcal{A}}_{i j}=0$, for $i \neq j \in\{1,2, \ldots, \ell \}.$
\end{cor}

\noindent
\textbf{Declaration:}  We declare that no conflicts of interest are associated with this publication and there has been no significant financial support for this work.  We certify that the submission is original work. 


\begin{thebibliography}{99}



\bibitem{Alekseevsky-Marchiafava} D.V. Alekseevsky and S. Marchiafava: 
Almost complex submanifolds of quaternionic manifolds, Steps in 
Differential Geometry, Institutes of Mathematics and Informatics, 
University of Debrecen, (2001), 23-38.

\bibitem{Apostol} T.M. Apostol: Mathematical analysis, second edition, 
Addison-Wesley Publishing Co., Reading, Mass.-London-Don Mills, Ont., 1974. xvii+492 pp.

\bibitem{Aydin} M. E. Aydin, A. Mihai, I. Mihai: Some Inequalities on submanifolds in statistical manifolds of constant curvature, Filomat, 29 (2015), no. 3, 465–477.

\bibitem{Ayi}  H. Aytimur and C. \"Ozg\"ur: Sharp inequalities for anti-invariant Riemannian submersions from Sasakian space forms, Journal of
 Geometry and Physics, 166 (2021), 104251.
 
\bibitem{Baird} P. Baird and J.C. Wood: Harmonic morphisms between 
Riemannian manifolds, London Mathematical Society Monographs, Vol. 29, Oxford University Press, The Clarendon Press, Oxford, (2003).

\bibitem{Besse} A.L. Besse: Einstein Manifolds, Ergebnisse der Mathematik 
und ihrer Grenzgebiete, 3, Folge 10, Springer, Berlin, (1987).


 \bibitem{Chen} B.-Y. Chen: Some pinching and classification theorems for minimal submanifolds, Arch. Math. (Basel), 60 (1993), 568–578.
 
\bibitem{Chen99} B.-Y. Chen: Relations between Ricci curvature and shape operator for submanifolds with arbitrary codimensions, Glasg. Math. Jour., 41 (1999), no. 1, 33–41.

 \bibitem{Chenbook} B.-Y. Chen: Pseudo-Riemannian Geometry, $\delta$-invariants and Applications, World Scientific Publishing Co. Pt.Ltd., Hackensack, NJ, (2011).

\bibitem{Eco} R.H. Escobales Jr: Riemannian submersion with totally geodesic
fibres, Jour. Differ. Geom., 10 (1967), 715-737.

\bibitem{Falcitelli-Ianus and Pastore} M. Falcitelli, S. Ianus and A.M.Pastore: Riemannian Submersions and related topics, World Scientific, 
River Edge, (2004).


\bibitem{Gulbahar} M. G\"ulbahar, \c S. E. Meri\c c and E. Kili\c c: Sharp inequalities involving the Ricci curvature for Riemannian submersions, Kragujevac ¨
Journal of Mathematics, 41 (2017), no. 2, 279–293.


\bibitem{Ianus-Mazzocco-Vilcu} S. Ianus, R. Mazzocco and G.E. Vilcu: Riemannian submersions from quaternionic manifolds, Acta. Appl. Math., 104 (2008), 83-89. 

\bibitem{S Ishi} S. Ishihara: Quaternion K\"ahlerian manifolds, J. Differ. Geom., 9 (1974), 483-500. 


\bibitem{Naghi}  M. Naghi, M. S. Stankovic and F. Alghamdi: Chen's improved inequality for pointwise hemi-slant warped products in K\"aehler manifolds, Filomat, 34 (2020), no. 3, 807-814.

\bibitem{O'Neill} B. O'Neill: The Fundamental Equations of a Submersion, 
Mich. Math. Jour., 13 (1966), 459-469.



\bibitem{Park-17} K.-S. Park: H-anti invariant submersions from almost 
quaternionic Hermitian manifolds, Czech. Math. Jour., 67 (2017), 557-558.



\bibitem{Sahin-10} B. \c Sahin: Anti-invariant Riemannian submersions from almost 
Hermitian manifolds, Cent. Eur. J. Math., 8 (2010), 437-447.

\bibitem{Sahin16}  B. \c Sahin: Chen’s first inequality for Riemannian maps, Ann. Polon. Math., 117 (2016), no. 3, 249–258.
\bibitem{MMT}  M. M. Tripathi: Chen-Ricci inequality for curvature like tensor and its applications, Diff. geom. Appl., 29 (2011), no. 5, 685–692.
 
\bibitem{Watson} B. Watson: Almost Hermitian submersions. Jour. Diff. Geom.,
11 (1976), no. 1, 147-165.


\end{thebibliography}
\end{document}